\documentclass[times]{acsauth}
\usepackage{amsmath,amsfonts,amssymb,amsthm,mathtools}

\newcommand{\tr}[1]{\mathop{{\rm \bf tr}\left[#1\right]}\nolimits}
\newcommand{\E}[1]{\mathop{{\rm \bf E}\left\{#1\right\}}\nolimits}
\theoremstyle{plain} %definition remark
\newtheorem{Algorithm}{Algorithm}
\newtheorem{remark}{Remark}
\newtheorem{lemma}{Lemma}
\newtheorem{example}{Example}

\begin{document}

\runningheads{M.~V.~Kulikova, J.~V.~Tsyganova}{Constructing stable adaptive Kalman filter-based techniques}

\title{Constructing numerically stable Kalman filter-based algorithms for gradient-based adaptive filtering}

\author{M.~V.~Kulikova\affil{1}\corrauth\,
J.~V.~Tsyganova\affil{2}}

\address{\affilnum{1}CEMAT, Instituto Superior T\'{e}cnico, Universidade de Lisboa, Portugal\break
\affilnum{2}Department of Mathematics and Information Technologies, Ulyanovsk State University, Russian Federation}

\corraddr{CEMAT, Instituto Superior T\'{e}cnico, Universidade de Lisboa, Av. Rovisco Pais 1,  1049-001, Lisbon, Portugal. E-mail: maria.kulikova@ist.utl.pt}

\begin{abstract}
This paper addresses the numerical aspects of adaptive filtering (AF) techniques for simultaneous state and parameters estimation arising in the design of dynamic positioning systems in many areas of research. The AF schemes consist of a recursive optimization procedure to identify the uncertain system parameters by minimizing an appropriate defined performance index and the application of the Kalman filter (KF) for dynamic positioning purpose. The use of gradient-based optimization methods in the AF computational schemes yields to a set of the filter sensitivity equations and a set of matrix Riccati-type sensitivity equations. The filter sensitivities evaluation is usually done by the conventional KF, which is known to be numerically unstable, and its derivatives with respect to unknown system parameters. Recently, a novel square-root approach for the gradient-based AF by the method of the maximum likelihood has been proposed. In this paper, we show that various square-root AF schemes can be derived from only two main theoretical results.  This elegant and simple computational technique replaces the standard methodology based on direct differentiation of the conventional KF equations (with their inherent numerical instability) by advanced square-root filters (and its derivatives as well). As a result, it improves the robustness of the computations against roundoff errors and leads to accurate variants of the gradient-based AFs. Additionally, such methods are ideal for simultaneous state estimation and parameter identification since all values are computed in parallel. The numerical experiments are given.
\end{abstract}


\maketitle

\vspace{-6pt}

\section{Introduction}
\vspace{-2pt}
The problem of developing the adaptive filtering (AF) techniques for simultaneous state and parameters estimation arising in the design of dynamic positioning systems has received increasing attention in recent years. Any AF method consists of a recursive optimization procedure to identify the uncertain system parameters by minimizing an appropriate defined performance index (e.g. the negative likelihood function) and the application of the Kalman filter (KF) for a dynamic positioning purpose. The gradient-based AF techniques additionally require the performance index (PI) gradient evaluation. It yields to a set of the filter sensitivity equations and a set of matrix Riccati-type sensitivity equations~\cite{Gupta1974,Mehra1974}. The sensitivities evaluation is usually done by the conventional KF and the direct differentiation of its equations (with respect to unknown system parameters); see~\cite{sandell1978maximum,Segal1988,segal1989new,hassani2013,Leander2014} and many others. A serious limitation of this methodology is the numerical instability of the conventional KF (with respect to round off errors) that may destroy the filter and, hence, the PI evaluation with the entire AF computational scheme.

Since 1960s there has been a great practical interest in the design of numerically stable and computationally efficient KF implementation methods. This has resulted in a large number of square-root (SR) filters, UD-based KF implementations and the fast Chandrasekhar-Kailath-Morf-Sidhu KF algorithms~\cite{Dyer1969,KaminskiBryson1971,Morf1974,Bierman1977,Sayed1994,ParkKailath1995}. Any of these advanced KF methods can replace the conventional KF in the AF schemes for a more stable PI evaluation. We may remark that current implementations of the KF are most often expressed in (what is called) an array square-root (ASR) form. They imply utilization of numerically stable orthogonal transformations for each recursion step. This feature enables more efficient parallel implementation and leads to algorithms with better numerical stability and conditioning properties; see \cite[Chapter~12]{KailathSayed2000} for an extended explanation.

Despite the existing diversity of the efficient KF algorithms, the PI \emph{gradient} evaluation (with respect to unknown system parameters) in terms of advanced KF methods is seldom addressed. In this paper we design simple and elegant computational scheme that allows for a natural extension of any ASR KF on the case of the filter sensitivities evaluation. Such methods are ideal for simultaneous state estimation and parameter identification since all values are computed in parallel. Additionally, our approach avoids implementation of the conventional KF (and its derivatives) because of its inherent numerical instability and, hence, improves the robustness of the computations. The first paper on a stable  filter sensitivity computation has suggested an extension of the information-type KF~\cite{Bierman1990}. Then,  the stable methods in terms of the covariance-type ASR KFs have been investigated in~\cite{2009Kul-IEEE,2009Kul-Matcom,2012TsyKul-AiT}. In this paper, we show that all types of the gradient-based AF schemes within stable ASR-based filters can be derived from two main theoretical results proven here. In contrast to the earlier published works, we do not derive a particular PI gradient evaluation method, but present a general approach that is able to extend any ASR KF (existing or new) on the robust filter derivatives computation. Additionally, the lower triangular scheme for the PI gradient evaluation is designed. This case has never been studied before. The numerical experiments are also given.

\section{State and parameter estimation of state-space models} \label{sec:MLE}
\vspace{-2pt}
Consider discrete-time linear stochastic system of the form
  \begin{align}
    x_{k} = & F(\theta) x_{k-1}+B(\theta)u_{k-1} + G(\theta) w_{k-1}, \quad  w_{k-1} \sim   {\mathcal N}(0, Q(\theta)),  \label{eq2.1} \\
    z_k   = & H(\theta) x_k+v_k,                   \quad  v_{k} \sim  {\mathcal N}(0, R(\theta))  \label{eq2.2}
  \end{align}
where $k$ is a discrete time ($k=1, \ldots, N$), i.e. $x_k$ means $x(t_k)$; vectors $x_k \in \mathbb R^n$ and $z_k \in \mathbb R^m$ are, respectively, the unknown dynamic state and the available measurements; $u_k \in \mathbb R^d$ is the deterministic input signal.  The  process noise, $\{w_k\}$, and  the measurement noise, $\{v_k\}$, are uncorrelated Gaussian white-noise processes, with covariance matrices $Q(\theta) \ge 0$ and $R(\theta) > 0$, respectively. All random variables
have known mean values, which we can take without loss of generality to be zero. The initial state $x_0$ is Gaussian random vector with the mean $\bar x_0(\theta)$ and the covariance matrix $\Pi_0(\theta)$, i.e. $x_0 \sim {\mathcal N}(\bar x_0(\theta), \Pi_0(\theta))$. It is independent from $\{w_k\}$ and $\{v_k\}$. Additionally, system~\eqref{eq2.1}, \eqref{eq2.2} is parameterized by a vector of unknown system parameters $\theta \in \mathbb R^p$, which needs to be estimated.  This means that the state-space model is known up to certain parameters, i.e. the matrices  $F(\theta) \in {\mathbb R}^{n\times n}$, $B(\theta) \in {\mathbb R}^{n\times d}$, $G(\theta) \in {\mathbb R}^{n\times q}$, $Q(\theta) \in {\mathbb R}^{q\times q}$, $H(\theta) \in {\mathbb R}^{m\times n}$ and $R(\theta) \in {\mathbb R}^{m\times m}$ may all depends on $\theta$. We stress that the initials conditions, i.e. $\bar x_0(\theta)$ and $\Pi_0(\theta) \in {\mathbb R}^{n\times n}$ may also depend on the parameters, however, such situation is seldom studied in the literature.

If there is no uncertainties in the system (i.e. $\theta$ is known and, hence, the state-space model is time-invariant), then the KF can be used for estimating the unobservable dynamic state $x_{k}$ from the corrupted measurements $z_1, \ldots, z_k$ as follows~\cite{KailathSayed2000}:
\begin{align}
\hat x_{k+1|k} & =  F \hat x_{k|k-1}+Bu_k + K_{p,k}e_k, & \hat x_{0|-1} & =  \bar x_0,  \label{state} && \\
 K_{p,k}       & = FP_{k|k-1}H^TR_{e,k}^{-1},           & e_k \; \;    & = z_k-H\hat x_{k|k-1},   &  R_{e,k}      & = HP_{k|k-1}H^T+R  \label{riccati1}
\end{align}
where $K_{p,k}=\E{\hat x_{k+1|k} e_k^T}$ and $e_k \sim {\cal N}(0, R_{e,k})$ are innovations of the discrete-time KF. The matrix $P_{k|k-1}$ appearing in the above formulas is the  error covariance matrix, i.e.
$P_{k|k-1}=\E{ (x_{k}-\hat
x_{k|k-1})(x_{k}-\hat x_{k|k-1})^{T}}$, and satisfies the difference
Riccati equation
\begin{align}
\label{riccati2}
 P_{k+1|k}  = & FP_{k|k-1}F^T+GQG^T - K_{p,k}R_{e,k}K_{p,k}^T, &  P_{0|-1} = & \Pi_0 > 0.
\end{align}
%Although the matrices $F$, $B$, $G$, $H$, $Q$ and $R$ are constant, the KF is still time-varying. It can be shown that
%observability is sufficient for convergence of the KF:
%if the pair $(F, H)$ is observable, then the Riccati recursion for $P_{k+1|k}$ converges to steady-state value $\lim \limits_{k \to \infty} P_{k+1|k} = %P_{\infty}$; see~\cite{Moore2005,Simon2006} for a detailed explanation.  We may also remark that if $P_{k+1|k}$ converges, then so do $K_{p,k}$ and $R_{e,k}$, %i.e. we get a stationary error distribution: $\lim \limits_{k \to \infty} e_k \to {\cal N } (0, R_{e,\infty})$ with finite variance $\lim \limits_{k \to \infty} %R_{e,k} = R_{e,\infty}$. Hence, the uncertainty in the state estimates remains bounded for all states. Next, if the state-space model is reachable, then the KF %reaches a steady-state value. Thus, under the assumptions that the pair $(F, H)$ is observable and the pair $(F, GQ^{1/2})$\footnote{The matrix $Q^{1/2}$ is the %square-root of $Q$, i.e. $Q=Q^{1/2}Q^{T/2}$.}   is reachable, the KF converges to the unique value and this stationary filter is asymptotically optimal. These %assumptions are also standard for most subspace identification algorithms~\cite{Jansson1998,Gustafsson2002}. 
In the next section we consider the problem of parameters estimation by the gradient-based AF techniques.

\subsection{Gradient-based adaptive filtering schemes}

The state-space model~\eqref{eq2.1}, \eqref{eq2.2} under examination is known up to certain parameters, $\theta \in \mathbb R^p$. This means that the associated KF~\eqref{state}~-- \eqref{riccati2} depends on the unknown $\theta$ as well. We stress that both the dynamic state, $x_k$, and system parameters, $\theta$, must be estimated simultaneously from only the observed noisy signal $z_k$. The classical way of solving such a problem  is to use \emph{adaptive} KF techniques, where the model parameters are estimated together with the dynamic state~\cite{Sarkka2009}.

To start implementing any AF scheme, one should choose first a PI that reflects the difference between the actual system and the utilized model with associated KF, which needs to be tuned up~\cite{hassani2013}. Then, a particular AF method is to be applied. At present, there are available many commonly used ways for the AF design in practice. Among them are the output-error techniques, the least-squares approach, the maximum-likelihood method, min-max entropy algorithms, {\it etc}~\cite{mehra1972approaches}. An important problem arising in this setting is convergence conditions of the constructed AF, i.e. convergence properties of the unknown parameter estimates, for both linear and nonlinear systems; see a consistency-oriented discussion in~\cite{Luders1974,Ljung1978,Bastin1988,marino1992global} and many others. For instance, \cite[Lemma~3.1]{Ljung1978} proves the main convergence result on this issue. It applies to quite a general situation and can be used as a common framework for the convergence and consistency analysis of many above-cited AF design methods. Throughout the paper we assume that all the assumptions of Lemma 3.1 hold; see details in~\cite[p.~776]{Ljung1978}.

The method of maximum likelihood is a general method for parameter estimation and often used in practice; see, for instance,~\cite{Gupta1974,Mehra1974,Segal1988,Bierman1990} and many others. It requires the maximization of the likelihood function (LF) given as follows~\cite{Schweppe1965}:
\begin{equation}
{\mathcal L}_{\theta}\left(Z_1^N\right)  =  -\frac{Nm}{2}\ln(2\pi) -
\frac{1}{2} \sum \limits_{k=1}^N \left\{
 \ln\left(\det R_{e,k}\right) + e_k^T  R_{e,k}^{-1}e_k \right\}  \label{LLF}
\end{equation}
where $Z_1^N=\{z_1,\ldots, z_N\}$ is $N$-step measurement history and $e_k \sim {\cal
N}\left(0,R_{e,k}\right)$ are the innovations generated by the discrete-time KF~\eqref{state}~-- \eqref{riccati2}.

Hence, the negative log LF represents the PI for solving the parameters estimation problem by the method of maximum likelihood. Then, a recursive  optimization procedure is used to identify the unknown system parameters $\theta$ by minimizing the PI. The optimization is often done by gradient-based or Newton's type methods where the computation of the LF gradient (LG) is necessary.  The basic iteration in gradient-type non-linear programming methods has the following form~\cite{Gupta1974}:
\begin{equation}
\label{optimization}
\theta_{n} = \theta_{n-1} - \gamma \left.\nabla \mu(\theta) \right|_{\theta = \theta_{n-1}}, \; n=1,2, \ldots
\end{equation}
where $\theta_n$ is the parameter vector at the $n$-th iteration and  $\left.\nabla \mu(\theta) \right|_{\theta = \theta_{n-1}}$ is the gradient of the PI with  respect to $\theta$ evaluated at $\theta = \theta_{n-1}$. The $\gamma$ is a scalar step size
parameter chosen to ensure that $\left.\mu(\theta) \right|_{\theta_{n}} \le \left. \mu(\theta) \right|_{ \theta_{n-1}} + \epsilon$ where $\epsilon$ is a positive number that can be chosen in a variety of ways; see~\cite{Gupta1974} for more details.

As can be seen, the gradient-based AF approach requires the run of the KF at each iteration step of the optimization method (i.e. for each $\theta_{n-1}$) to generate the $\{ e_k, R_{e,k} \}$, $k = 1, \ldots, N$ for the PI evaluation, $\left.\mu(\theta)\right|_{\theta = \theta_{n-1}}$, corresponding to the current approximation $\theta_{n-1}$. Additionally, it demands the gradient computation, $\left.\nabla \mu(\theta)\right|_{\theta = \theta_{n-1}}$ at each $\theta_{n-1}$. This leads to a set of $p$ vector equations, known as the {\it filter sensitivity equations}, and a set of $p$ matrix equations, known as the {\it Riccati-type sensitivity equations}. The described {\it forward filter method} demands roughly an implementation of  $p+1$ equivalent KF's all running in the forward time direction where $p$ is a number of the unknown system parameters.

In this manuscript, we do not discuss the particular optimization method that can be applied in each particular situation, but explain how the PI (the negative log LF) and its gradient can be computed accurately together with the system state. Such methods are ideal for simultaneous state estimation and parameter identification since all values are calculated in parallel.

\subsection{The problem of numerically instability of the conventional KF}

Both parts of the AF scheme, i.e. the chosen optimization method (for finding the optimal $\hat \theta^*$) and the chosen KF algorithm (for computing the PI and estimating $x_k$), play an important role in the computational scheme and affect the accuracy of the recursive adaptive estimator. Most of the previously proposed AF techniques are based on the conventional KF~\eqref{state}~-- \eqref{riccati2} and the direct differentiation of its equations for the PI gradient evaluation~\cite{sandell1978maximum,Segal1988,segal1989new,hassani2013,Leander2014}. The main disadvantage of this approach is numerical instability of the conventional KF while the requirement to compute the filter sensitivities in parallel deteriorates the situation. Here, we improve the accuracy of gradient-based AF methodology by replacing the numerically unstable conventional KF to advanced KF methods and their derivatives with respect to unknown system parameters. More precisely, we are focusing in the techniques developed in the KF community to solve ill conditioned problems. To start the presentation of our main results, we first discuss the ASR filters.

The matrix $P_{k|k-1}$ appearing in~\eqref{state}~-- \eqref{riccati2} has the physical meaning of being the variance of the state prediction error, $x_k-\hat x_{k|k-1}$, and therefore has to be nonnegative-definite. Round off errors may destroy this property leading to a failure of the filter. In contrast to the conventional KF~\eqref{state}~-- \eqref{riccati2},  the ASR methods propagate only square-root factors\footnote{Throughout the paper we use the Cholesky decomposition of the form $A=A^{T/2} A^{1/2}$, where $A^{1/2}$ is an upper triangular matrix  with positive diagonal entries.
 %We also introduce the notation $A^{-1/2}=(A^{1/2})^{-1}$, $A^{T/2}=(A^{1/2})^T$ and $A^{-T/2}=(A^{-1/2})^T$.
 }
  $P_{k|k-1}^{1/2}$ of the covariance matrices $P_{k|k-1}$, $k=1, \ldots, N$. The point is that the product of the computed factors, say  $\hat P_{k|k-1}=\hat P_{k|k-1}^{T/2}\hat P_{k|k-1}^{1/2}$, is a symmetric matrix with positive elements on the diagonal and it is almost certainly nonnegative-definite; see~\cite[Chapter~12]{KailathSayed2000} for more details. Furthermore, any ASR filter uses a numerically stable orthogonal rotation at each iteration step. This feature enables more efficient parallel implementation and leads to algorithms with better numerical stability and conditioning properties.

All types of the ASR implementations can be divided into two simple cases. Some of them uses the orthogonal transformation of the form $QA=R$ with $R$ being an upper triangular matrix and others imply the transformation $QA=L$ where $L$ is a lower triangular matrix\footnote{ The left-hand side matrix $A$ is called the pre-array of the ASR filter. The right-hand side matrices $R$ and $L$ are called the post-arrays.}. We illustrate this statement by two ASR KF algorithms designed in~\cite{ParkKailath1995}.

\textsc{The extended square-root covariance filter} (\textbf{eSRCF}). Given the initial values for the filter: $P_{0|-1}^{-T/2} \hat x_{0|-1}=\Pi_{0}^{-T/2} \bar x_{0} $ and $P_{0|-1}^{1/2}=\Pi_0^{1/2}$, recursively update ($k=1, \ldots, N$):
\begin{align}
Q
\left[
\begin{array}{cc|c}
R^{1/2} & 0 & -R^{-T/2}z_k \\
P_{k|k-1}^{1/2}H^T & P_{k|k-1}^{1/2}F^T &
P_{k|k-1}^{-T/2} \hat x_{k|k-1} \\
0 &  Q^{1/2}G^T &  0
\end{array}
\right]
 & =
 \left[
 \begin{array}{cc|c}
R_{e,k}^{1/2} & \bar K_{p,k}^T  & -\bar e_k \\
0 & P_{k+1|k}^{1/2} &  P_{k+1|k}^{-T/2} \hat x_{k+1|k} \\
0 & 0 &  \gamma_k
\end{array}
\right],\label{eq:esrcf:1}\\
\hat x_{k+1|k} & =  \left( P_{k+1|k}^{T/2}\right) \left(
P_{k+1|k}^{-T/2}\hat x_{k+1|k}\right) \label{eq:esrcf:2}
\end{align}
where $Q$ is any orthogonal transformation  such that the first two (block) columns of the matrix on the right-hand side of formula~\eqref{eq:esrcf:1} is upper
triangular. We introduce a notation for the normalized innovations $\bar e_k = R_{e,k}^{-T/2}e_k$ and the normalized Kalman gain $\bar K_{p,k} = FP_{k|k-1}H^TR_{e,k}^{-1/2}$. The matrix $R_{e,k}^{1/2}$ is a square-root factor of $R_{e,k}$.

\begin{remark}
\label{remark:1}
The parentheses in~\eqref{eq:esrcf:2} are used to indicate the quantities
that can be directly read off from the post-array in~\eqref{eq:esrcf:1}. Hence, no matrices need to be inverting for finding the state vector estimate $\hat x_{k+1|k}$, $k=1, \ldots, N$.
\end{remark}

\textsc{The extended square-root information filter} (\textbf{eSRIF}). Given the initial values for the filter: $P_{0|-1}^{-T/2} \hat x_{0|-1}=\Pi_{0}^{-T/2} \bar x_{0} $ and $P_{0|-1}^{-T/2}=\Pi_0^{-T/2}$, recursively update ($k=1, \ldots, N$):
\begin{eqnarray}
& Q &
\left[
\begin{array}{ccc|c}
R^{-T/2} & -R^{-T/2}HF^{-1} & R^{-T/2}HF^{-1}GQ^{T/2}
& -R^{-T/2}z_k  \\
0 & P_{k|k-1}^{-T/2}F^{-1}  & -P_{k|k-1}^{-T/2}F^{-1}GQ^{T/2}
& P_{k|k-1}^{-T/2}\hat x_{k|k-1} \ \\
0 & 0 & I  & 0
\end{array}
\right]  \nonumber \\
& = & \left[
\begin{array}{ccc|c}
R_{e,k}^{-T/2} & 0 & 0  & -\bar e_k  \\
-P_{k+1|k}^{-T/2} K_{p,k} &
P_{k+1|k}^{-T/2} &  0  &P_{k+1|k}^{-T/2} \hat x_{k+1|k}  \\
* & * &  *  & *
\end{array}
\right] \label{eq:esrif:1}
\end{eqnarray}
where $Q$  is any orthogonal transformation  such that the first
three (block) columns of the post-array is a lower triangular matrix. The predicted  state estimate can be
found by solving the triangular system of the following form:
\begin{equation}
\left( P_{k+1|k}^{-T/2}\right) \left(
\hat x_{k+1|k}\right)= \left(P_{k+1|k}^{-T/2}\hat x_{k+1|k}\right). \label{eq:esrif:2}
\end{equation}

\begin{remark}
\label{remark:2}
The eSRCF and eSRIF  can be verified by ``squaring'' both sides of the $QA=R$ (or $QA=L$), using the fact that
$QQ^T=I$, and comparing the entries of both sides of the result.
The detailed derivations can be also found in~\cite{KailathSayed2000}.
\end{remark}

As mentioned earlier, the maximum likelihood estimation procedure leads to implementation of the KF (and its derivatives with respect to unknown system parameters), which is known to be numerically unstable. It is desirable to avoid the use of the conventional KF in the computational scheme. In other words, we would like to replace the disadvantageous conventional KF by numerically stable ASR filters, e.g. by the eSRCF/eSRIF presented above. The log LF and its gradient can be expressed in terms of the quantities appearing in the ASR filters as follows~\cite{2009Kul-IEEE}:
\begin{align}
{\mathcal L}_{\theta}\left(Z_1^N\right)  = & -\frac{Nm}{2}\ln(2\pi) -
\frac{1}{2} \sum \limits_{k=1}^N \left\{2
 \ln\left(\det R_{e,k}^{1/2}\right) + \bar e_k^T \bar e_k \right\},  \label{LLF:sqrt} \\
\frac{\partial {\mathcal L}_{\theta}\left(Z_1^N
\right)}{\partial \theta_i}  = & -\sum \limits_{k=1}^N  \left\{
  \tr{ R_{e,k}^{-1/2} \cdot
\frac{\partial  R_{e,k}^{1/2}}{\partial \theta_i}}  +   \bar e_k^T \frac{ \partial \bar e_k}{\partial \theta_i}  \right\},
\quad i=1, \ldots, p   \label{grad-esrcf}
\end{align}
where ${\rm \bf tr}[ \cdot ]$ denotes the trace of matrices.

In the next section, we design a simple and convenient technique for computing derivatives of the ASR filter variables required in equation~\eqref{grad-esrcf}.

\section{ASR Filter Derivatives Computation} \label{sec:main_result}
\vspace{-2pt}
 First, we note that each iteration of ASR filters has the following form:  $QA=B$ where $Q$ is any orthogonal transformation such that  the post-array $B$ is either a lower triangular  or upper triangular matrix. Treating these two cases separately, we prove the following main results.

\begin{lemma} [\textsc{The lower triangular case}]
\label{lemma-1}
 Let entries of the pre-array $A \in {\mathbb
R}^{(s+k)\times (s+l)}$ be known differentiable
functions of a parameter $\theta$. Consider the equation of the form $QA=L$ with the following partitioning:
\begin{equation}
\label{assume-2}  \begin{array}{rcc} &
\begin{array}{cc}
s & \; \; \; \; l
 \end{array}  & \\
 Q &  \left[
 \begin{array}{c|c}
 A_{11} & A_{12} \\
 A_{21} & A_{22}
 \end{array}
 \right]
 & \!\!\!\!\!\!\!\!
\begin{array}{cc}
 k \\
 s
 \end{array}
\end{array}
 \begin{array}{rcc}
 &  \begin{array}{cc}
 s & \; \; \; \; l
 \end{array} & \\
 = & \left[
 \begin{array}{c|c}
 0 & L_{12}\\
 L_{21} & L_{22}
 \end{array}
 \right] &
 \!\!\!\!\!\!\!\!
\begin{array}{c}
 k \\
 s
 \end{array}
\end{array}
\end{equation}
where $Q \in {\mathbb R}^{(s+k)\times (s+k)}$ is an orthogonal
matrix that lower-triangularizes the first (block) column of the matrix on the left-hand side
of~(\ref{assume-2}) and $L_{21} \in {\mathbb R}^{s\times s}$ is
lower triangular. Introduce the notation
\begin{equation}
\label{notation-lemma:1}     \begin{array}{rcc} &
  \begin{array}{cc}
 s & \; \; \;\; l
 \end{array}  & \\
 Q &  \left[
 \begin{array}{c|c}
\left( A_{11}\right)'_{\theta} & \left( A_{12}\right)'_{\theta} \\
 \left( A_{21}\right)'_{\theta} & \left( A_{22}\right)'_{\theta}
 \end{array}
 \right]
 & \!\!\!\!\!\!\!\!
  \begin{array}{c}
 k \\
 s
 \end{array}
\end{array}
   \begin{array}{rcc}
&   \begin{array}{cc}
 s & \; \; \; \;  l
 \end{array} & \\
 = & \left[
 \begin{array}{c|c}
 X & N \\
 Y & V
 \end{array}
 \right] &
 \!\!\!\!\!\!\!\!
 \begin{array}{c}
 k \\
 s
 \end{array}
\end{array}\;.
\end{equation}
 Then  given the derivative of the pre-array $A'_{\theta}$, the following formulas calculate the corresponding derivatives of the post-array blocks:
\begin{equation}
 \label{lemma2:eq:3}
 \left(L_{21}\right)'_{\theta}=(\bar
{\cal U}^T + {\cal D} + \bar {\cal L})L_{21},
\end{equation}
\begin{equation}
\label{lemma2:eq:4} \left(L_{22}\right)'_{\theta}  =  \left[ \bar {\cal U}^T -\bar {\cal U} \right] L_{22} +
L_{21}^{-T}X^T L_{12} +   V
\end{equation}
where  $\bar {\cal L}$, ${\cal D}$ and $\bar {\cal U}$ are respectively strictly lower
triangular, diagonal and strictly upper triangular parts of the following matrix product $YL_{21}^{-1}$.
\end{lemma}

\begin{proof}
At first, we show that $Q'_{\theta}Q^T$ is a skew symmetric matrix.
For that, we differentiate both sides of the formula $QQ^T=I$ with
respect to $\theta$ and arrive at
$Q'_{\theta}Q^T+Q\left(Q^T\right)'_{\theta}=0$, or in the
equivalent form $Q'_{\theta}Q^T=-\left(Q'_{\theta}Q^T\right)^T$.
The latter implies that the matrix $Q'_{\theta}Q^T$ is skew
symmetric and can be presented as a difference of two matrices, i.e.
$Q'_{\theta}Q^T = \bar {\cal U}^T- \bar {\cal U}$ where $\bar {\cal U}$ is an $(s+k)\times (s+k)$ strictly
upper triangular matrix. Thus, the last $(s\times s)$-block located at the main diagonal of
$Q'_{\theta}Q^T$ has the same form, i.e.
\begin{equation}
\label{new:new}
\left[Q'_{\theta}Q^T\right]_{s\times s} = \bar {\cal U}^T_{s\times s}- \bar {\cal U}_{s\times s}
\end{equation}
 where $\bar {\cal U}_{s\times s}$ is a $s \times s$ strictly upper triangular matrix and $\left[Q'_{\theta}Q^T\right]_{s\times s}$ stands for the $(s\times
s)$-matrix composed of the entries located at the intersections of
the last $s$ rows with the last $s$ columns of the product $Q'_{\theta} Q^T$.

Next, we prove that the above-mentioned matrix $\bar {\cal U}_{s\times s}$  is, in fact, the upper triangular part of the
matrix product $YL_{21}^{-1}$. To do this, we
differentiate the first equation in formula~(\ref{assume-2}), i.e.
$$
 \begin{array}{rcc}
 &  \begin{array}{c}
 s
 \end{array}  & \\
 Q &  \left[
 \begin{array}{c}
 A_{11}  \\
 A_{21}
 \end{array}
 \right]
 & \!\!
 \begin{array}{c}
 k \\
 s
 \end{array}
\end{array}
\begin{array}{rcc}
 &  \begin{array}{c}
 s
 \end{array} & \\
 = & \left[
\begin{array}{c}
 0  \\
 L_{21}
 \end{array}
 \right] &
 \!\!
\begin{array}{c}
 k \\
 s
 \end{array}
\end{array}
$$
with respect to $\theta$. Then, taking into account
notation~(\ref{notation-lemma:1}) and equality $A=Q^TL$, we obtain
\begin{equation}
\label{proof2:1} \left[ \begin{array}{c}
 0 \\
 \left(L_{21}\right)'_{\theta}
 \end{array}
 \right]
 =  Q'_{\theta} \left[
\begin{array}{c}
 A_{11}  \\
 A_{21}
 \end{array}
 \right]
+Q \left[
\begin{array}{c}
\left( A_{11}\right)'_{\theta}  \\
 \left(A_{21}\right)'_{\theta}
 \end{array}
 \right]
=Q'_{\theta}Q^T \left[ \begin{array}{c}
 0  \\
 L_{21}
 \end{array}
 \right]
+\left[\begin{array}{c}
X  \\
Y
 \end{array}\right].
 \end{equation}

Further, it is not difficult to see that the pseudoinverse matrix (Moore-Penrose inversion)  of
$\left[\; 0 \; | \; L_{21}\right]^T$  is $\left[ \; 0 \; | \; L_{21}^{-1}\right]$. Therefore the right multiplication of both sides
of~(\ref{proof2:1}) by the pseudoinverse yields
\begin{equation}
\label{proof2:2} \left[
 \begin{array}{c|c}
 0 & 0\\
\hline
 0 & \left(L_{21}\right)'_{\theta} L^{-1}_{21}
 \end{array}
 \right]=
 Q'_{\theta}Q^T\left[0_{k\times k} \oplus I_{s\times s} \right]+\!
\left[\begin{array}{c}
X  \\
 Y
 \end{array}\right]\!
\left[ 0 \; L_{21}^{-1}\right]\!
\end{equation}
where $I_{s\times s}$ is the identity matrix of dimension~$s$ and
$0_{k\times k}$ is the zero block of size $k \times k$. The $\left[0_{k\times k} \oplus I_{s\times s}\right]$
means $\mbox{\rm diag}\{0_{k\times k},I_{s\times s}\}$. Now we
remark that
\begin{equation}
\label{proof2:3}
 \left[
 \begin{array}{c|c}
0 & 0\\[3pt]
\hline
 0 &  \left(L_{21}\right)'_{\theta} L^{-1}_{21}
 \end{array}
 \right]  =
 \left[
 \begin{array}{c|c}
0 & \left[Q'_{\theta}Q^T\right]_{col:~last~s}^{row:~first~k}\\[3pt]
 \hline
  0 &  \left[Q'_{\theta}Q^T\right]_{s\times s}
 \end{array}
 \right] + \left[
 \begin{array}{c|c}
0 & XL_{21}^{-1} \\[3pt]
 \hline
 0 & YL_{21}^{-1}
 \end{array}
 \right]
\end{equation}
where  $\left[Q'_{\theta}Q^T\right]_{col:~last~s}^{row:~first~k}$
stands for the $(k\times s)$-matrix composed of the entries located
at the intersection of the first $k$ rows with the last $s$ columns
of the matrix $Q'_{\theta}Q^T$.

From the matrix equation~(\ref{proof2:3}), we conclude the
following. First, the matrix on the left-hand side
of~(\ref{proof2:3}) is block lower triangular. Thus, the strictly
upper triangular part of the matrix
$\left[Q'_{\theta}Q^T\right]_{s\times s}$ must exactly annihilate
the strictly upper triangular part of the corresponding second term
on the right-hand side of~(\ref{proof2:3}). In other words, if the
matrix product $YL_{21}^{-1}$ is represented as
$$
 YL_{21}^{-1}=\bar {\cal L}_{s\times s}+{\cal D}_{s\times s}+\bar {\cal U}_{s\times s}
$$
where $\bar {\cal L}_{s\times s}$, ${\cal D}_{s\times s}$ and $\bar {\cal U}_{s\times s}$ are respectively the strictly
lower triangular, diagonal and strictly upper triangular parts, then
the matrix $\bar {\cal U}_{s\times s}$, in fact, satisfies~(\ref{new:new}).

Now formula~(\ref{lemma2:eq:3}) is easily justified. Indeed, from
the matrix equation~(\ref{proof2:3}), we obtain
\begin{eqnarray}
 \left(L_{21}\right)'_{\theta} L^{-1}_{21} & = & \underbrace{\bar {\cal U}^T_{s\times s} - \bar
{\cal U}_{s\times s}}_{\left[Q'_{\theta}Q^T\right]_{s\times s}} + \underbrace{\bar
{\cal L}_{s\times s} + {\cal D}_{s\times s} + \bar {\cal U}_{s\times s}}_{ YL_{21}^{-1}}, \nonumber \\
\left(L_{21}\right)'_{\theta} & = & (\bar
{\cal U}^T_{s\times s} + {\cal D}_{s\times s} + \bar {\cal L}_{s\times s})L_{21}.
\nonumber
\end{eqnarray}
For the sake of simplicity, in equation~(\ref{lemma2:eq:3})  we omit the subscripts of the matrices $\bar {\cal L}$, ${\cal D}$ and $\bar {\cal U}$.

Second, from the matrix equation~(\ref{proof2:3}) we observe that the first (block) row of the left-hand side matrix
in~(\ref{proof2:3}) is zero. Thus, the first (block) row of the
matrix $Q'_{\theta}Q^T$ must exactly cancel the corresponding block of
the second term in~(\ref{proof2:3}), i.e. we arrive at
\begin{equation}
 \label{part_new:1}
\left[Q'_{\theta}Q^T\right]_{col:~last~s}^{row:~first~k} = -XL_{21}^{-1}.
\end{equation}

Next, we wish to validate~(\ref{lemma2:eq:4}). By differentiating the
last equation in~(\ref{assume-2}) with respect to $\theta$, and then taking into account
notation~(\ref{notation-lemma:1}), we derive
$$
 \left[
\begin{array}{c}
\left(L_{12}\right)'_{\theta} \\
\left(L_{22}\right)'_{\theta}
\end{array}
\right]  =  Q'_{\theta} \left[
\begin{array}{c}
A_{11} \\
A_{21}
 \end{array}
\right]
  + Q \left[
\begin{array}{c}
\left(A_{12}\right)'_{\theta} \\
\left(A_{22}\right)'_{\theta}
 \end{array}
\right]= Q'_{\theta}Q^T \left[
\begin{array}{c}
L_{12} \\
L_{22}
 \end{array}
\right]
  +  \left[
\begin{array}{c}
N \\
V
 \end{array}
\right].
$$

The previous formula implies that
\begin{eqnarray}
\label{th:proof:3} \left(L_{22}\right)'_{\theta} &  =  & V +  \left[ Q'_{\theta}Q^T
\right]_{col:~first~k}^{row:~last~s} L_{12}+\left[
Q'_{\theta}Q^T \right]_{s\times s}L_{22} \nonumber \\
& =  & V +  \left[ Q'_{\theta}Q^T
\right]_{col:~first~k}^{row:~last~s} L_{12} +  \left[ \bar {\cal U}_{s\times s}^T -\bar {\cal  U}_{s\times s}  \right] L_{22}
\end{eqnarray}
where $\bar {\cal  U}_{s\times s}$ is the upper triangular matrix from~(\ref{new:new}) and $\left[ Q'_{\theta} Q^T
\right]_{col:~first~k}^{row:~last~s}$ stands for the $(s\times
k)$-matrix composed of the entries located at the intersections of
the last $s$ rows with the first $k$ columns of the product
$Q'_{\theta} Q^T$.

Eventually, formula~(\ref{part_new:1}) and the fact that
$Q'_{\theta} Q^T$ is skew symmetric result in
  \begin{equation}
  \label{th:proof:1}\!\!
 \left[Q'_{\theta} Q^T \right]_{col:~first~k}^{row:~last~s}  =
  - \left[\!\left[ Q'_{\theta} Q^T \!\right]_{col:~last~s}^{row:~first~k}\right]^T\! \!\!= - \left[ -
 XL_{21}^{-1}\right]^T =  L_{21}^{-T}X^T.
  \end{equation}

Thus, the substitution of~(\ref{th:proof:1}) in~(\ref{th:proof:3})
validates~(\ref{lemma2:eq:4}) and completes the proof of
Lemma~\ref{lemma-1}.
\end{proof}

\begin{lemma} [\textsc{The upper triangular case}]
\label{lemma-2} Let entries of the pre-array $A \in {\mathbb
R}^{(s+k)\times (s+l)}$ be known differentiable
functions of a parameter $\theta$. Consider the equation of the form $QA=R$ with the following partitioning:
\begin{equation}
\label{assume-2-2}  \begin{array}{rcc} &
\begin{array}{cc}
s & \; \; \; \; l
 \end{array}  & \\
 Q &  \left[
 \begin{array}{c|c}
 A_{11} & A_{12} \\
 A_{21} & A_{22}
 \end{array}
 \right]
 & \!\!\!\!\!\!
\begin{array}{cc}
 s \\
 k
 \end{array}
\end{array}
 \begin{array}{rcc}
 &  \begin{array}{cc}
 s & \; \; \; \; l
 \end{array} & \\
 = & \left[
 \begin{array}{c|c}
 R_{11} & R_{12} \\
 0 & R_{22}
 \end{array}
 \right] &
 \!\!\!\!\!\!\!\!
\begin{array}{c}
 s \\
 k
 \end{array}
\end{array}
\end{equation}
where $Q \in {\mathbb R}^{(s+k)\times (s+k)}$ is an orthogonal
matrix that produces the block zero entry on the right-hand side
of~(\ref{assume-2-2}) and $R_{11} \in {\mathbb R}^{s\times s}$ is
upper triangular. Introduce the notation~(\ref{notation-lemma:1}). Then  given the derivative of the pre-array $A'_{\theta}$, the following formulas calculate the corresponding derivatives of the post-array:
\begin{equation}
 \label{lemma2-1}
 \left(R_{11}\right)'_{\theta}=(\bar
{\cal L}^T + {\cal  D}  + \bar {\cal  U})R_{11},
\end{equation}
\begin{equation}
\label{lemma2-2} \left(R_{12}\right)'_{\theta}  =   \left[ \bar
{\cal  L}^T -\bar {\cal  L} \right] R_{12} +
R_{11}^{-T}Y^T R_{22}+N
\end{equation}
where  $\bar {\cal  L}$, ${\cal  D}$ and $\bar {\cal U}$ are respectively strictly lower
triangular, diagonal and strictly upper triangular parts of the following matrix product $XR_{11}^{-1}$.
\end{lemma}

\begin{proof}
Lemma~\ref{lemma-2} can be proved at the same way as Lemma~\ref{lemma-1}. The detail derivation of the formulas above can be also found in~\cite{2009Kul-Matcom}.
\end{proof}

\section{Summary of the Computations}
\vspace{-2pt}

Theoretical results presented in Lemmas~\ref{lemma-1}, \ref{lemma-2} yield a general computational scheme for the filter derivative computations. This new approach is able to replace the conventional KF (and its derivatives with respect to unknown system parameters) by any numerically stable ASR filter in the gradient-based AF techniques. The ASR methodology utilizes the pre-array $A$ of any chosen ASR filter and its derivatives in order to compute the post-array and its derivatives, respectively. Algorithms~\ref{alg:1}, \ref{alg:2} summarize the entire computational schemes in details.

\noindent
\begin{tabular}{p{0.5\textwidth}|p{0.5\textwidth}}
\begin{Algorithm} \label{alg:1} {\small (\textsc{The lower triangular case}) }

\vspace{0.3cm}
\underline{Input Data:} The pre-array $A$ and its derivatives $\partial A/ \partial \theta_i$, $i=1, \ldots, p$.

\vspace{0.2cm}
 \underline{Process:} Compute the post-array $L$ by~(\ref{assume-2}). Save matrices $\{ Q, L\}$ for future steps. Then, for each component $\theta_i$, $i=1, \ldots, p$:
\begin{itemize}
\item Find $Q\displaystyle\frac{\partial A}{\partial \theta_i}$ and introduce the notations as in~(\ref{notation-lemma:1}). Save the blocks  $\{ X_i, Y_i, N_i, V_i\}$;
\item Calculate $Y_iL_{21}^{-1}$. Split it into strictly lower triangular $\bar {\cal L}_i$, diagonal ${\cal D}_i$ and strictly upper triangular $\bar {\cal U}_i$ parts;
\item Compute
$\displaystyle\frac{\partial
L_{21}}{\partial \theta_i}=\bigl(\bar {\cal U}_i^T+{\cal D}_i+\bar {\cal L}_i\bigr)L_{21};$
\item $\displaystyle\frac{\partial
L_{22}}{\partial \theta_i} \! =\!  \left[ \bar {\cal U}_i^T \!\! -\!\bar {\cal U}_i \right] L_{22}\! +\!
L_{21}^{-T}X_i^T L_{12}\! + \! V_i$.
\end{itemize}
\underline{Output Data:} The post-array $L$ and its derivatives: $\partial L_{21}/ \partial \theta_i$, $\partial L_{22}/ \partial \theta_i$, $i=1, \ldots, p$.
\end{Algorithm}
&

\begin{Algorithm} \label{alg:2} {\small (\textsc{The upper triangular case})}

\vspace{0.3cm}
\underline{Input Data:} The pre-array $A$
 and its derivatives $\partial A/ \partial \theta_i$, $i=1, \ldots, p$.

\vspace{0.2cm}
\underline{Process:} Compute the post-array $R$ by~(\ref{assume-2-2}). Save matrices $\{ Q, R\}$ for future steps. Then, for each component $\theta_i$, $i=1, \ldots, p$:
\begin{itemize}
\item Find $Q\displaystyle\frac{\partial A}{\partial \theta_i}$ and introduce the notations as in~(\ref{notation-lemma:1}).  Save the blocks  $\{ X_i, Y_i, N_i, V_i\}$;
\item Calculate $X_iR_{11}^{-1}$. Split it into strictly lower triangular $\bar {\cal L}_i$, diagonal ${\cal D}_i$ and strictly upper triangular $\bar {\cal U}_i$ parts;
\item Compute $\displaystyle\frac{\partial
R_{11}}{\partial \theta_i}=\bigl(\bar {\cal L}^T_i+{\cal D}_i+\bar {\cal U}_i\bigr)R_{11};$
\item $\displaystyle\frac{\partial
R_{12}}{\partial \theta_i} \! = \!\left[ \bar {\cal L}_i^T \!\! - \! \bar {\cal L}_i \right]
R_{12} \! + \! R_{11}^{-T}Y_i^TR_{22}\!+ \!N_i$.
\end{itemize}
  \underline{Output Data:} The post-array $R$ and its derivatives: $\partial R_{11}/ \partial \theta_i$ and $\partial R_{12}/ \partial \theta_i$, $i=1, \ldots, p$.
\end{Algorithm}
\end{tabular}
\vspace{0.3cm}

Having applied Algorithms~\ref{alg:1}, \ref{alg:2} at each iteration of the ASR KF, we obtain the post array of the filter and its derivatives with respect to unknown system parameters for each $k=1, \ldots, N$. These quantities contain the $\{ e_k, R_{e,k} \}$ and $\{ \partial e_k/\partial \theta_i, \partial R_{e,k}/\partial \theta_i \}$, $i=1, \ldots p$, required for the PI and its gradient evaluation; see~\eqref{LLF:sqrt}, \eqref{grad-esrcf}. Hence, the entire gradient-based AF computational scheme can be formulated as follows. Let  $\theta_{n-1}$ denotes the value of $\theta$ after $n-1$ iterations of the optimization algorithm~\eqref{optimization}. In this section we explain how the next cycle for computing $\theta_n$ can be obtained by using the chosen gradient-based optimization method, the chosen PI and any ASR filter, e.g. the eSRCF/eSRIF presented above.

\begin{Algorithm} \label{alg:3} (\textsc{Adaptive filtering scheme})

\vspace{0.3cm}
\noindent
\underline{Input Data:} A current approximation $\theta_{n-1}$.

\vspace{0.2cm}
\noindent
\underline{Process:} Evaluate the system matrices (and its derivatives) at the current $\theta_{n-1}$: $\hat F(\theta)=F(\theta)\left|_{\theta_{n-1}}\right.$,
$\hat G(\theta) =G(\theta)\left|_{\theta_{n-1}} \right.$ {\it etc}. To improve robustness of the computations, replace the unstable conventional KF~(\ref{state})~-- (\ref{riccati2}) by any ASR filtering algorithm. Use the Cholesky decomposition to find the square-root of the matrices: $\hat \Pi_{0}^{1/2}$ and $\hat R^{1/2}$, $\hat Q^{1/2}$. Set the initial values for the filter and, then, process the measurements $\{ z_1, \ldots, z_N\}$ as follows:
\begin{itemize}
\item Form the pre-array and its derivatives of the chosen ASR filter.
\item Given the pre-array (and its derivatives), find the post-array and its derivatives (with respect to each $\theta_i$, $i=1, \ldots, p$) as follows.
%%%%%%%%%%%%%%%%%%%%%%%%%%%%%%%%%%%%%%%%%%%%%%%%%%%55
%If the pre-array is a lower triangular matrix, then apply Algorithm~\ref{alg:1}. If the pre-array is an upper triangular matrix, then apply Algorithm~\ref{alg:2}.
If the post-array has the form of a lower triangular matrix, then apply Algorithm~\ref{alg:1}.
If the post-array has the form of an upper triangular matrix, then apply Algorithm~\ref{alg:2}.
%%%%%%%%%%%%%%%%%%%%%%%%%%%%%%%%%%%%%%%%%%%%%%%%
\item Extract $\bar e_k$ and $R_{e,k}^{1/2}$ from the post-array. Compute new term in the PI.
\item Extract $\partial \bar e_k/ \partial \theta_i$ and
$\partial R_{e,k}^{1/2} / \partial \theta_i $, $i=1, \ldots, p$ from the derivatives of the post-array. Compute new term in the PI gradient.
\end{itemize}
 After processing all measurements $\{ z_1, \ldots, z_N\}$, the PI and its gradient are evaluated. Next, use the chosen gradient-based method in order to find the next approximation $\theta_n$.

\vspace{0.2cm}
\noindent
\underline{Output Data:} Next approximation $\theta_{n}$.
\end{Algorithm}

Repeat Algorithm~\ref{alg:3} for the next $\theta_{n+1}$ ($n=1, 2, \ldots$) until the stopping criterion is satisfied. The proposed technique simultaneously identifies the uncertain system parameters by minimizing the PI and estimates the unknown state vector of dynamic system.

\section{Illustrative examples: the eSRCF- and eSRIF-based AF methods}\label{sec:method}
\vspace{-2pt}
The detailed derivation of the eSRCF-based technique for the log LF and its gradient evaluation can be found in~\cite{2009Kul-IEEE}. Here we show how the method can be easily obtained from Lemma~\ref{lemma-2} and Algorithms~\ref{alg:2}. First, we note that the post-array of the eSRCF filter is an upper triangular matrix. Next, the matrix that needs to be triangularized is of size $n+m$. Hence, we apply Lemma~\ref{lemma-2} to the eSRCF pre-array with  $s=m+n$, $k=q$, $l=1$ and the following partitioning:
$$
Q\;
\underbrace{
\begin{bmatrix}
\boxed{\begin{matrix}
R^{1/2} & 0 \\
P_{k|k-1}^{1/2}H^T & P_{k|k-1}^{1/2} F^T
\end{matrix}}
&
\boxed{\begin{matrix}
 -R^{-{\rm T}/2}z_k \\
P_{k|k-1}^{-{\rm T}/2} \hat x_{k|k-1}
\end{matrix}} \\
 \mbox{\tiny $A_{11} \in {\mathbb R}^{(m+n)\times (m+n)}$}  &  \mbox{\tiny $A_{12} \in {\mathbb R}^{(m+n)\times 1}$} \\
\boxed{\begin{matrix}
\qquad 0 \qquad & Q^{1/2}G^{\rm T} \; \;\;
\end{matrix}}
 &
\boxed{\begin{matrix}
\qquad \; 0 \; \qquad \;\; \;
\end{matrix}}\\
 \mbox{\tiny $A_{21} \in {\mathbb R}^{q\times (m+n)}$}  &  \mbox{\tiny $A_{22} \in {\mathbb R}^{q\times 1}$}
\end{bmatrix}
}_{Pre-array \; A}
=
\underbrace{
\begin{bmatrix}
\boxed{ \begin{matrix}
R_{e,k}^{1/2} & \bar K_{p,k}^T   \\
0 & P_{k+1|k}^{1/2}
\end{matrix}} &
\boxed{\begin{matrix}
-\bar e_k \\
P_{k+1|k}^{-T/2} \hat x_{k+1|k}
\end{matrix}} \\
 \mbox{\tiny $R_{11} \in {\mathbb R}^{(m+n)\times (m+n)}$}  &  \mbox{\tiny $R_{12} \in {\mathbb R}^{(m+n)\times 1}$} \\
\boxed{\begin{matrix}
\quad 0 \quad & \quad 0 \; \; \;\;
\end{matrix}}
 &
\boxed{ \begin{matrix} \qquad \gamma_k \qquad \; \; \;
\end{matrix}} \\
 \mbox{\tiny $R_{21} \in {\mathbb R}^{q\times (m+n)}$}  &  \mbox{\tiny $R_{22} \in {\mathbb R}^{q\times 1}$}
\end{bmatrix}
}_{Post-array \; R}.
$$
The computational scheme of Algorithm~\ref{alg:2} leads to the filter derivative computations and, in particular, to the  $\partial R_{e,k}^{1/2}/ \partial \theta_i$ and $\partial \bar e_{k}/ \partial \theta_i$, $i=1, \ldots, p$  evaluation required in the PI and its gradient evaluation.

\vspace{0.3cm}
At the same way the information-type algorithm can be easily obtained from the eSRIF; see also the detailed derivation for the log LF and its gradient evaluation  in~\cite{2006KulSem-LN}. We note that the post-array of the eSRIF filter is a lower triangular matrix. Hence, we apply Lemma~\ref{lemma-1} and Algorithm~\ref{alg:1} to the eSRIF with $s=m+n+q$, $k=0$, $l=1$ and the following partitioning:
\begin{align*}
Q &
\underbrace{
\begin{bmatrix}
\boxed{\begin{matrix}
\phantom{R^{-T/2}} & \phantom{-R^{-T/2}HF^{-1}} & \phantom{R^{-T/2}HF^{-1}GQ^{T/2} }
\end{matrix}}
&
\boxed{\begin{matrix}
\phantom{P_{k|k-1}^{-T/2}\hat x_{k|k-1}}
\end{matrix}} \\
 \mbox{\tiny $A_{11}$ is empty}  &  \mbox{\tiny $A_{12}$ is empty} \\
\boxed{\begin{matrix}
R^{-T/2} & -R^{-T/2}HF^{-1} & R^{-T/2}HF^{-1}GQ^{T/2}  \\
0 & P_{k|k-1}^{-T/2}F^{-1}  & -P_{k|k-1}^{-T/2}F^{-1}GQ^{T/2}\\
0 & 0 & I
\end{matrix}}
 &
\boxed{\begin{matrix}
-R^{-T/2}z_k     \\
P_{k|k-1}^{-T/2}\hat x_{k|k-1} \\
0
\end{matrix}}\\
 \mbox{\tiny $A_{21} \in {\mathbb R}^{(m+n+q)\times (m+n+q)}$}  &  \mbox{\tiny $A_{22} \in {\mathbb R}^{(m+n+q)\times 1}$}
\end{bmatrix}
}_{Pre-array \; A} \\
= &
\underbrace{
\begin{bmatrix}
\boxed{ \begin{matrix}
\phantom{-P_{k+1|k}^{-T/2} K_{p,k}} & \phantom{P_{k+1|k}^{-T/2}} & \phantom{ 0}
\end{matrix}} &
\boxed{\begin{matrix}
\phantom{P_{k+1|k}^{-T/2} \hat x_{k+1|k} }
\end{matrix}} \\
 \mbox{\tiny $L_{11}$ is empty}  &  \mbox{\tiny $L_{12}$ is empty} \\
\boxed{\begin{matrix}
R_{e,k}^{-T/2} & 0 & 0    \\
-P_{k+1|k}^{-T/2} K_{p,k} & P_{k+1|k}^{-T/2} &  0   \\
* & * &  *
\end{matrix}}
 &
\boxed{ \begin{matrix}
-\bar e_k    \\
P_{k+1|k}^{-T/2} \hat x_{k+1|k} \\
*
\end{matrix}} \\
 \mbox{\tiny $L_{21} \in {\mathbb R}^{(m+n+q)\times (m+n+q)}$}  &  \mbox{\tiny $L_{22} \in {\mathbb R}^{(m+n+q)\times 1}$}
\end{bmatrix}
}_{Post-array \; L}.
\end{align*}

In summary, the proposed computational schemes naturally extend any ASR filter and allow {\it the filter and the filter sensitivity equations} to be updated in parallel. Hence, such methods are ideal for simultaneous state estimation and parameter identification.

\begin{remark}
Some modern ASR KF implementations are based on the $UDU^T$ factorization of the pre-array. Hence, an alternative approach to a problem of numerically stable PI and its gradient evaluation can be found in, the so-called, UD-based filters developed first in~\cite{Bierman1977}. The problem of the UD-based filters' derivative computation (with respect to unknown system parameters) has been formulated by Bierman {\it et al.} in~\cite{Bierman1990} and has been open since 1990s. It was recently solved in~\cite{Tsyganova2013IEEE}.
\end{remark}

\section{Numerical Examples\label{experiments}}
\vspace{-2pt}

First, we wish to check our theoretical derivations presented in Lemma~\ref{lemma-1} and~\ref{lemma-2}. To do so,
we consider the following simple test problems.

\begin{example}\label{ex:1:1} (\textsc{Simple test problem: the upper triangular case})

\vspace{0.3cm}
\noindent
For the given pre-array
$$
 A= \left[
\begin{array}{ccc|c}
{\theta^5}/20 \; & \; \theta^4/8 & {\theta^3}/{6} & \;  \theta^3/3\\
{\theta^4}/8 \; & \; \theta^3/3 & {\theta^2}/{2} & \;  {\theta^2}/{2}\\
{\theta^3}/{6} \; & \; {\theta^2}/{2} & \theta & \; 1
\end{array}
\right],
$$
compute the post-arrays $R$ and its derivative $R'_{\theta}$, say, at $\theta=2$ where the first three (block) columns of the post-array $R$ is an upper triangular matrix.
\end{example}

We note, that the unknown parameter $\theta$ is a scalar value, i.e. $p=1$. For simplicity, we assume that $N=1$, i.e. we illustrate the
detailed explanation of only one iteration step of the algorithm. Next, we remark that the post-array should be an upper triangular matrix and, hence, Lemma~\ref{lemma-2} and Algorithm~\ref{alg:2} should be applied to solve the stated problem. Then, we pay an attention to the partitioning in~(\ref{assume-2-2}) from  Lemma~\ref{lemma-2} and conclude that $s=3$, $l=1$ and $k=0$. Hence, the blocks $A_{21}$, $A_{22}$ of the pre-array $A$ and, respectively, the $R_{21}$, $R_{22}$ of the post-array $R$ are empty. Indeed, according to Example~\ref{ex:1:1} the first three (block) columns of the post-array $R$ is an upper triangular matrix. This means that $s=3$ and, hence, $k=0$, i.e. $A_{21}$, $A_{22}$ are empty. As a result, $l=1$.

Having applied the computational scheme from Algorithm~\ref{alg:2} to the pre-array in Example~\ref{ex:1:1}, we compute the post-array $R$ and its derivative (at the point $\theta=2$). The obtained results are summarized in Table~\ref{tab:1}. All codes were written in MATLAB. To check our derivations, we compute the norm $\left|\left|{(A^T A)}'_{\theta =2}-{(R^T R)}'_{\theta=2}\right|\right|_{\infty}$. Indeed, from equation $QA=R$ we have $A^T A = R^TR$. Thus, the derivatives of both sides of the latter formula must also agree. The obtained value is $1.33\cdot 10^{-14}$. This confirms the correctness of the calculation of Algorithm~\ref{alg:2} and validates the theoretical derivations of Lemma~\ref{lemma-2}.

\begin{table}
\caption{Numerical results for Example~\ref{ex:1:1}}  \label{tab:1}
\centering
\tabsize
\begin{tabular}{|ll|}
\toprule
\multicolumn{2}{|l|}{  We are given the pre-array $A$ and its derivatives with respect to each $\theta_i$, (in the example $p=1$):} \\
 & Pre-array
$
 A= \left[
\begin{array}{ccc|c}
{\theta^5}/20 \; & \; \theta^4/8 & {\theta^3}/{6} & \;  \theta^3/3\\
{\theta^4}/8 \; & \; \theta^3/3 & {\theta^2}/{2} & \;  {\theta^2}/{2}\\
{\theta^3}/{6} \; & \; {\theta^2}/{2} & \theta & \; 1
\end{array}
\right]$, i.e. $\left.A\right|_{\theta=2}=\left[
\begin{array}{rrr|r}
 1.6000  &  2.0000 &   1.3333 &   2.6667 \\
    2.0000  &  2.6667  &  2.0000 &   2.0000 \\
    1.3333  &  2.0000  &  2.0000 &   1.0000
\end{array}
\right]$, \\
 & and
$A'_{\theta}=\left[
\begin{array}{ccc|c}
\theta^4/4 & \theta^3/2 & \theta^2/2 & \theta^2 \\
\theta^3/2 & \theta^2   & \theta     & \theta \\
\theta^2/2 & \theta     &    1      &     0
\end{array}
\right]$. So,  $\left.A'_{\theta}\right|_{\theta=2}=\left[
\begin{array}{ccc|c}
     4   &  4    & 2    & 4 \\
     4   &  4    & 2    & 2 \\
     2   &  2    & 1    & 0
\end{array}
\right]$. \\
\midrule
\multicolumn{2}{|l|}{Compute the post-array $R$ using $QR$ algorithm and save matrices $\{ Q, R\}$ for future steps:} \\
 & Post-array
$R=\left[
\begin{array}{rrr|r}
 -2.8875   & -3.8788 &   -3.0476 &   -3.3247\\
         0 &   -0.2576 &   -0.6954 &    0.8886\\
         0  &        0  &   0.0797  &   0.5179
\end{array}
\right]$,
$Q=
\begin{bmatrix}
  -0.5541  & -0.6926 &  -0.4618 \\
    0.5795 &   0.0773 &   -0.8113 \\
    0.5976 &  -0.7171 &    0.3586
 \end{bmatrix}$. \\
\midrule
\multicolumn{2}{|l|}{Apply the designed derivative computation method ($p=1$):} \\
$\bullet$  & Compute
$QA'_{\theta}$. Denote $X_1=\begin{bmatrix}
 -5.9105  & -5.9105  & -2.9552 \\
    1.0045 &   1.0045 &   0.5022 \\
    0.2390 &   0.2390 &   0.1195
\end{bmatrix}$, $\begin{matrix}
Y_1 = [\quad] \\
V_1 = [\quad] \\
\end{matrix}
$,
$N_1 =
\begin{bmatrix}
 -3.6017 \\
    2.4725 \\
    0.9562
\end{bmatrix}$.\\
$\bullet$  & Find
$X_1R_{11}^{-1}= \begin{bmatrix}
  2.0469  & -7.8778 & -27.5511 \\
   -0.3479 &   1.3388 &    4.6822 \\
   -0.0828 &   0.3186 &    1.1143
\end{bmatrix}$. Split it into $\bar {\cal L}_1 =
\begin{bmatrix}
  0  & 0 & 0 \\
   -0.3479 & 0 &    0 \\
   -0.0828 &   0.3186 & 0
\end{bmatrix}$, \\
 & \phantom{Find} ${\cal D}_1= \begin{bmatrix}
  2.0469  & 0 & 0 \\
  0 &   1.3388 &   0 \\
  0 &   0 &    1.1143
\end{bmatrix}$, $\bar {\cal U}_1 =
\begin{bmatrix}
  0  & -7.8778 & -27.5511 \\
   0 &  0 &  4.6822 \\
   0 &  0 & 0
\end{bmatrix}$.\\
$\bullet$  & Calculate $\left.R'_{11}\right|_{\theta=2}=
\begin{bmatrix}
 -5.9105 &  -5.8209 &  -2.7199   \\
         0 &   -0.3448 &  -0.5325 \\
         0 &        0  &  0.0888
\end{bmatrix}$ and $\left.R'_{12}\right|_{\theta=2}=
\begin{bmatrix}
 -3.9537 \\
 1.4810 \\
 0.3978
\end{bmatrix}$. \\
\midrule
\multicolumn{2}{|l|}{Hence, the derivative of the post-array is} \\
 & $\left.R'_{\theta}\right|_{\theta=2}= \left[
\begin{array}{rrr|r}
 -5.9105 &  -5.8209 &  -2.7199  & -3.9537 \\
         0 &   -0.3448 &  -0.5325 &    1.4810 \\
         0 &        0  &  0.0888  &  0.3978
\end{array}
\right]$. \\
\bottomrule
\multicolumn{2}{|l|}{Accuracy of  the computations: $\left|\left|{(A^T A)}'_{\theta
=2}-{(R^T R)}'_{\theta=2}\right|\right|_{\infty}=1.33\cdot 10^{-14}$} \\
\bottomrule
\end{tabular}
\end{table}

\begin{example}\label{ex:1:2} (\textsc{Simple test problem: the lower triangular case})

\vspace{0.3cm}
\noindent
 For the pre-array $A$ from example~\ref{ex:1:1}, compute the post-arrays $L$ and its derivative $L'_{\theta}$ (at $\theta=2$) where the first three (block) columns of the post-array $L$ is a lower triangular matrix; see equation~(\ref{assume-2}).
\end{example}

The lower triangular case can be justified at the same way. We note, that $l=1$,  $s=3$, $k=0$ and, hence, we have the partitioning~(\ref{assume-2}) of the pre-array $A$ with the empty blocks $A_{11}$, $A_{12}$. The post-array $L$ is block lower triangular and, hence, we apply  the computational scheme presented in Algorithm~\ref{alg:1}. The obtained results are summarized in Table~\ref{tab:matrix:2}. The accuracy of the computation is $\left|\left|{(A^T A)}'_{\theta
=2}-{(L^T L)}'_{\theta=2}\right|\right|_{\infty}=2.57\cdot 10^{-14}$. This confirms the correctness of the calculation of Algorithm~\ref{alg:1} and validates the theoretical derivations of Lemma~\ref{lemma-1}.

\begin{table}
\caption{Numerical results for Example~\ref{ex:1:2}}  \label{tab:matrix:2}
\centering
\tabsize
\begin{tabular}{|ll|}
\toprule
\multicolumn{2}{|l|}{We are given the pre-array $A$ and its derivatives with respect to each $\theta_i$, (in the example $p=1$):} \\
 & Pre-array
$
 A= \left[
\begin{array}{ccc|c}
{\theta^5}/20 \; & \; \theta^4/8 & {\theta^3}/{6} & \;  \theta^3/3\\
{\theta^4}/8 \; & \; \theta^3/3 & {\theta^2}/{2} & \;  {\theta^2}/{2}\\
{\theta^3}/{6} \; & \; {\theta^2}/{2} & \theta & \; 1
\end{array}
\right]$, i.e. $\left.A\right|_{\theta=2}=\left[
\begin{array}{rrr|r}
 1.6000  &  2.0000 &   1.3333 &   2.6667 \\
    2.0000  &  2.6667  &  2.0000 &   2.0000 \\
    1.3333  &  2.0000  &  2.0000 &   1.0000
\end{array}
\right]$, \\
 & and
$A'_{\theta}=\left[
\begin{array}{ccc|c}
\theta^4/4 & \theta^3/2 & \theta^2/2 & \theta^2 \\
\theta^3/2 & \theta^2   & \theta     & \theta \\
\theta^2/2 & \theta     &    1      &     0
\end{array}
\right]$. So,  $\left.A'_{\theta}\right|_{\theta=2}=\left[
\begin{array}{ccc|c}
     4   &  4    & 2    & 4 \\
     4   &  4    & 2    & 2 \\
     2   &  2    & 1    & 0
\end{array}
\right]$.  \\
\midrule
\multicolumn{2}{|l|}{Compute the post-array $L$ using $QL$ algorithm and save matrices $\{ Q, L\}$ for future steps:} \\
 & Post-array
$L=\left[
\begin{array}{rrr|r}
   -0.0306   &      0   &      0  &  -0.6882 \\
   -0.6456   & -0.6195  &      0  & -1.5163  \\
   -2.8142   & -3.8376  & -3.1269 &  -3.0559
\end{array}
\right]$,
$Q=
\begin{bmatrix}
  -0.6882  &  -0.5869 &  -0.4264 \\
    0.6882 &  -0.3424 &  -0.6396 \\
   -0.2294 &   0.7337 &  -0.6396
 \end{bmatrix}$. \\
\midrule
\multicolumn{2}{|l|}{Apply the designed derivative computation method ($p=1$):} \\
$\bullet$  & Compute
$QA'_{\theta}$. Denote $\begin{array}{c}
X_1 = [\quad] \\
N_1 = [\quad] \\
\end{array}
$, $Y_1=
\begin{bmatrix}
 -0.4588   & -0.4588  &  -0.2294 \\
   -2.2499 &  -2.2499 &  -1.1250 \\
   -5.5432 &  -5.5432 &  -2.7716
\end{bmatrix}$,
$V_1 =
\begin{bmatrix}
   -1.3765 \\
   -3.0325 \\
   -2.9848 \\
\end{bmatrix}$.\\
$\bullet$  & Find
$X_1L_{21}^{-1}= \begin{bmatrix}
  2.2105   &  0.2861  &  0.0734 \\
   10.8396 &   1.4031 &   0.3598 \\
   26.7057 &   3.4569 &   0.8864
\end{bmatrix}$. Split it into $\bar {\cal L}_1 =
\begin{bmatrix}
  0  & 0 & 0 \\
   10.8396 & 0 &    0 \\
  26.7057 &   3.4569 &   0
\end{bmatrix}$, \\
 & \phantom{Find} ${\cal D}_1= \begin{bmatrix}
  2.2105 & 0 & 0 \\
  0 &   1.4031 &   0 \\
  0 &   0 &    0.8864
\end{bmatrix}$, $\bar {\cal U}_1 =
\begin{bmatrix}
  0  &  0.2861  &  0.0734 \\
   0 &  0 &   0.3598 \\
   0 &  0 & 0
\end{bmatrix}$.\\
$\bullet$  & Calculate $\left.L'_{21}\right|_{\theta=2}=
\begin{bmatrix}
  -0.0676  &       0  &       0 \\
   -1.2462 &  -0.8693 &        0 \\
   -5.7777 &  -5.7661 &  -2.7716
\end{bmatrix}$ and $\left.L'_{22}\right|_{\theta=2}=
\begin{bmatrix}
  -0.7184 \\
   -2.1301 \\
   -3.5808
\end{bmatrix}$. \\
\midrule
\multicolumn{2}{|l|}{Hence, the derivative of the post-array is} \\
 & $\left.L'_{\theta}\right|_{\theta=2}= \left[
\begin{array}{rrr|r}
   -0.0676  &       0  &      0  &  -0.7184 \\
   -1.2462  & -0.8693  &      0  & -2.1301 \\
   -5.7777  & -5.7661  & -2.7716 &  -3.5808
\end{array}
\right]$. \\
\bottomrule
\multicolumn{2}{|l|}{Accuracy of  the computations: $\left|\left|{(A^T A)}'_{\theta
=2}-{(L^T L)}'_{\theta=2}\right|\right|_{\infty}=   2.57\cdot 10^{-14}$} \\
\bottomrule
\end{tabular}
\end{table}

\vspace{0.3cm}
Next, we wish to discuss the convergence of the parameter $\theta$ to its real value, i.e. to discuss the accuracy of the designed recursive AF estimator presented in Algorithm~\ref{alg:3}. As mentioned earlier, the new AF scheme is developed from the techniques designed in the Kalman filtering community to solve ill conditioned problems. This should improve accuracy and robustness of the computations for a finite-precision computer arithmetics. To check this property, we consider the set of ill-conditioned test problems from~\cite{Tsyganova2013IEEE}.

\begin{example}\label{ex:3} (\textsc{Set of ill-conditioned test problems})

\vspace{0.3cm}
\noindent
Consider the state-space model~(\ref{eq2.1})-(\ref{eq2.2})
with $\{F, G, B, H, \Pi_0, Q, R \}$ given by
\begin{align*}
F= & \begin{bmatrix}
   1 & 0 & 0\\
   0 & 1 & 0 \\
   0 & 0 & 1
   \end{bmatrix},  B=\begin{bmatrix}
     0 \\
     0 \\
     0
     \end{bmatrix},
     G=\begin{bmatrix}
   0 \\
   0 \\
   0
   \end{bmatrix},
      Q=\begin{bmatrix}
   1
   \end{bmatrix}, \quad
          R= \begin{bmatrix}
                        \delta^2\theta^2 & 0  \\
                        0 & \delta^2\theta^2
                      \end{bmatrix},
                      H= \begin{bmatrix}
      1 & 1 &  1\\
      1 & 1 &  1+\delta
      \end{bmatrix} \\
      \mbox{ with } & x_0 \sim {\cal N}\left(\begin{bmatrix}
                                                        0 \\
                                                        0 \\
                                                        0
                                                        \end{bmatrix},
                                                        \begin{bmatrix}
                                                        \theta^2 & 0 & 0 \\
                                                        0 & \theta^2 & 0 \\
                                                        0 & 0 & \theta^2
                                                        \end{bmatrix}\right)
\end{align*}
where $\theta$ is an unknown system parameter, that needs to be
estimated. To
simulate roundoff we assume that $\delta^2<\epsilon_{\text{roundoff}}$, but
$\delta>\epsilon_{\text{roundoff}}$ where $\epsilon_{\text{roundoff}}$ denotes the
unit roundoff error\footnote{Computer roundoff for floating-point
arithmetic is often characterized by a single parameter
$\epsilon_{\text{roundoff}}$, defined in different sources as the largest
number such that either $1+\epsilon_{\text{roundoff}} = 1$ or
$1+\epsilon_{\text{roundoff}}/2 = 1$ in machine precision.}, i.e. the machine precision limit.
\end{example}

The set of ill-conditioned problems is constructed as follows. When $\theta  = 1$, Example~\ref{ex:3} coincides with well-known test from~\cite{GrewalAndrews2001} that demonstrates how a problem that is well conditioned, as posed, can be made ill-conditioned by the filter. It is often used in the Kalman filtering community for observing the influence of round off errors on various KF implementations. The difficulty is in matrix inversion $R_{e,k}$. After processing only the first measurement $z_1$, the matrix $ R_{e,1}=R+H\Pi_{0}H^T$ becomes singular in machine precision, i.e.  as  $\delta \to
\epsilon_{\text{roundoff}}$. This yields the failure of the conventional KF. To construct a proper test problem for the gradient-based AF estimators,  the authors of~\cite{Tsyganova2013IEEE}  introduced an unknown
system parameter $\theta$, making sure that the same problem is now applied to the matrix $(R_{e,1})'_{\theta}$. In other words,  for any fixed value of the parameter $\theta \ne 0$, the matrices $R_{e,1}=R+H\Pi_{0}H^T$ and $(R_{e,1})'_{\theta}$ are ill-conditioned in machine precision, i.e. as $\delta \to
\epsilon_{\text{roundoff}}$.  As a consequence, both parts of the gradient-based AF techniques (the PI and its gradient evaluation, respectively) fail after processing the first measurement.  This destroys the entire AF estimator grounded in the conventional KF implementation. Hence, such test allows for observing the influence of the round off errors on various gradient-based AF schemes.

%______________________________________________________________________
\begin{figure}
\includegraphics[width=\textwidth]{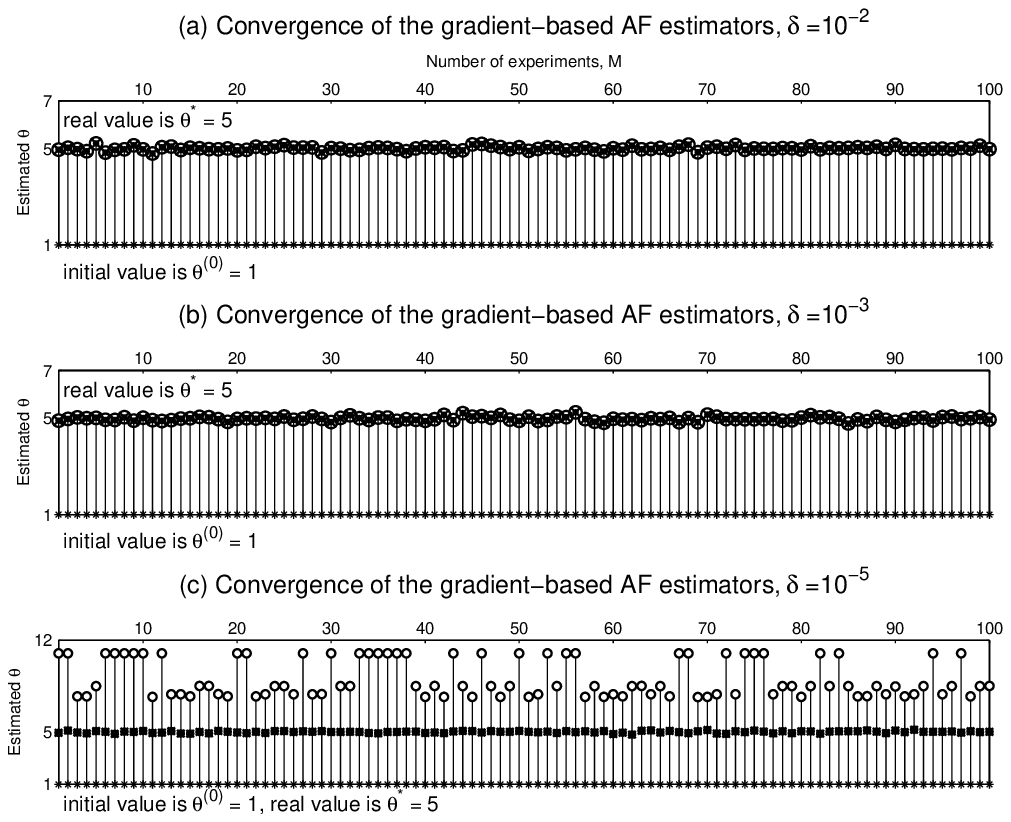}
\caption{The computed maximum likelihood estimates of $\theta$ by three gradient-based AF techniques: within conventional KF (marker $\circ$); the eSRCF implementation (marker $\bullet$) and the eSRIF filer (marker $\times$). The initial parameter value, i.e. $\theta^{(0)}=1$, is marked by $*$.} \label{fig:2}
\end{figure}
%______________________________________________________________________

We perform the following set of numerical experiments. Given the ``true'' value of the parameter
$\theta$, say $\theta^*=5$, the system is simulated for $1000$
samples for various values of $\delta$ while $\delta \to \epsilon_{\text{roundoff}}$. The generated
data is then used to solve the inverse problem, i.e. to compute the maximum likelihood estimates by gradient-based AF schemes. We consider the AF recursive estimator based on the conventional KF, on the eSRCF and eSRIF.  The designed  Algorithms~\ref{alg:1}, \ref{alg:2} are used for the PI and its gradient evaluation within numerically stable ASR filters (the eSRCF and eSRIF). Algorithm~\ref{alg:3} represents the general gradient-based AF scheme where we implemented the standard MATLAB built-in function \verb"fminunc" for optimization purpose. This optimization function utilizes the PI (the negative Log LF) and its gradient that are calculated by the conventional KF approach and the designed ASR methodology. The same initial value of $\theta^{(0)}=1$ is applied in all examined AF estimators. To observe the convergence of the parameter $\theta$ from the initial value $\theta^{(0)}=1$ to its real value $\theta^*=5$, we perform $100$ Monte Carlo simulations and illustrate the obtained results by Fig.~\ref{fig:2}.

From the first two graphs in Fig.~\ref{fig:2} we see that when $\delta=10^{-2}$ and $\delta=10^{-3}$, i.e. when the considered problem is well-posed, all gradient-based AF techniques work equally well.  We can observe their perfect convergence from the initial value $\theta^{(0)}=1$ to the real value $\theta^*=5$ in all $100$ Monte Carlo simulations. However, the situation dramatically changes for $\delta=10^{-5}$ when the problem becomes moderately ill-conditioned. The gradient-based AF scheme within conventional KF exhibits perfect performance for $\delta=10^{-2}$ and $\delta=10^{-3}$, but it completely fails for $\delta=10^{-5}$. Indeed, the conventional approach leads to incorrect parameter estimate in most cases among $100$ Monte Carlo simulations when $\delta=10^{-5}$. Meanwhile, the AF techniques based on the numerically stable ASR implementations work well for all examined $\delta$ as $\delta \to \epsilon_{\text{roundoff}}$.

\section{Conclusion}

In this paper, we developed an elegant and
simple general computational scheme that extends functionality of any array
square-root Kalman filtering algorithm on the filter derivative computations. These values are required in the gradient-based adaptive filtering  techniques for simultaneous state and parameter estimation of dynamic positioning systems in many areas of research. The proposed approach yields the improved robustness of the computations against roundoff errors.

\ack The first author thanks the support of Portuguese National Fund ({\it Funda\c{c}\~{a}o para a
Ci\^{e}ncia e a Tecnologia})  within the scope of project SFRH/BPD/64397/2009.
The authors also would like to express their gratitude to the anonymous referees for their valuable remarks and comments on the paper.

\bibliographystyle{wileyj}
\bibliography{BibTex_Library/list_books,%
              BibTex_Library/list_identification,%
              BibTex_Library/list_filters,%
              BibTex_Library/list_ML,%
              BibTex_Library/List-Kulikova}

\end{document}